\documentstyle[12pt]{article}
\begin{document}

\Large
\centerline{\bf Nonlocal Invariants in Index Theory}
\bigskip\bigskip

\large
\centerline{\bf Steven Rosenberg${}^\ast$}
\normalsize
\medskip
\centerline{\parbox{5in}{
 Department of Mathematics, Boston University,
Boston, MA 02215}} 

\centerline{sr@math.bu.edu}

\bigskip
\bigskip

\footnote{${}^\ast$Partially
supported by the NSF.}

\def\dvol{{\rm dvol}}

\normalsize

In its original form, the Atiyah-Singer Index Theorem equates
two global quantities of a closed manifold, one analytic 
(the index of an elliptic operator) and one topological
(a characteristic number).  Because it relates invariants from
different branches of mathematics,
the Index Theorem
has  many applications and extensions
to differential geometry, K-theory,
mathematical physics, and other fields.
This report focuses on advances in geometric
aspects of
index theory.

For operators naturally associated to a Riemannian metric on a closed
manifold, the topological side of the Index Theorem can
often be expressed as the
integral of local (i.e.~pointwise) curvature expression.  
We will first discuss
these local refinements in \S1, which arise naturally in heat equation
proofs of the Index Theorem.  In \S\S2,3, we
 discuss further developments in index theory which
lead to spectral invariants, the eta invariant and the 
determinant of an elliptic operator, that are definitely nonlocal.
Finally, in \S4 we point out some recent connections among these
nonlocal invariants and classical index theory.

\bigskip
\noindent \large{\bf \S1 Local invariants in Atiyah-Singer index
theory}

\bigskip
\normalsize
The Atiyah-Singer Index Theorem, first proved around 1963, equates the
index, \hfill\break
${\rm Ind} \ D = {\rm dim\ Ker }\ D - {\rm dim\ Coker}\ D$,
of an elliptic operator $D:\Gamma(E)\to\Gamma(F)$ taking
sections of a bundle $E$ 
over a
closed manifold $M$ to sections of a bundle $F$
with $\langle P(M,\sigma_{\rm top}(D)),[M]\rangle$,
a characteristic number built from the
topology of $M$ and topological information contained in the top order
symbol of $D$ \cite{a-s}, \cite{palais}.  For particular choices of
$D$, this deep result encompasses Chern's generalization of the
Gauss-Bonnet theorem, the Hirzebruch signature theorem, and
Hirzebruch's generalization of the Riemann-Roch theorem, as well as
giving many new results.   Atiyah and Singer also proved the Families
Index Theorem \cite{atiyah-singer;families} for a family of elliptic
operators
$\{D_n\}$ parametrized by $n$ in a compact manifold $N$.  This theorem
identifies  the Chern character of the index bundle ${\sc Ind}\ D$ 
 in $H^\ast (N;{\bf Q})$ with a characteristic class
on $N$ built from the topology of $N$ and the pushforward of the
symbols of the $D_n$; the significance of this theorem is that the
Chern character of a bundle determines the K-theory isomorphism class
of the bundle up to torsion.
  Here ${\sc Ind}\ D = {\sc Ker}\ D -
{\sc Coker }\ D$ is a virtual bundle whose
fiber for generic $n\in N$ is the formal difference ${\sc Ker}\ D_n -
{\sc Coker}\ D_n$.
For example, for a family of Dirac operators
associated to a family of metrics $g_n$ on $M$ parametrized by $n\in
N$,  we get ${\rm ch}({\sc Ind}\ D) = \pi_\ast \hat A(M/N)$, where
$\hat A(M/N)$ is the $\hat A$-polynomial of the bundle over $N$ whose
fiber
at $n$ is the bundle of spinors associated to $g_n,$   and
$\pi:M\times N\to N$ is the projection.  This example becomes less
trivial if we replace $M\times N$ by a manifold $X$ fibered by $M$
(cf.~\S4). 

Among the basic examples of  elliptic operators are the so-called
geometric operators, i.e. operators canonically associated
to a Riemannian metric on $M$.  For these operators, Hodge theory
identifies the index with the Euler characteristic of
 de Rham-type cohomology groups.
For example, the
Gauss-Bonnet operator is given by $D = d+ d^\ast$ taking even forms to
odd
forms, where $d$ is exterior differentiation and $d^\ast$ is its
adjoint with respect to the Hodge inner product on forms induced by
the metric.  The index of $D$ is the Euler characteristic $\chi(M)$.
The Dolbeault operator $\bar\partial + \bar\partial^\ast$,
 the analog of the Gauss-Bonnet operator for a
complex K\"ahler manifold, has the arithemetic genus as its index, and
the
signature operator has $\sigma(M)$, the signature of $M$, as its
index.  (The index of the Dirac operator on a spin manifold does not
have a Hodge theoretic topological interpretation.) 

The simplest geometric example is the Gauss-Bonnet operator on a
closed surface.
The characteristic number given by the Index Theorem 
is $\langle e(M), [M]\rangle$, the value
of the Euler class of $TM$ evaluated on the fundamental class of $M$.
By Chern-Weil theory, this characteristic number equals
$(1/2\pi)\int_M KdA$, the integral of the Gaussian curvature.  Thus
the Index Theorem reduces in this case to the classical Gauss-Bonnet
theorem: $\chi(M) = (1/2\pi) \int_M KdA.$   In particular, we have
expressed the index as the integral of a {\it local invariant}, a top
dimensional form (or density) computed pointwise canonically from the
metric
(i.e. independent of choice of chart).  

This derivation of the Gauss-Bonnet theorem is somewhat indirect, as
it uses Hodge theory to equate the analytic term, the index, with a
topological expression, and then  appeals to Chern-Weil theory to
obtain a geometric expression for the topological side of the Index
Theorem.
  Similar remarks apply to other
geometric operators.
Beginning with \cite{atiyah-bott-patodi},
heat equation
proofs of the Index Theorem were developed which expressed the index
of a geometrically defined operator directly as the
integral of the correct Chern-Weil expression.  
In brief,  it is easily shown that the supertrace ${\rm tr}_s
e^{-t{\bf D}^\ast {\bf D}}  = {\rm tr}(e^{-tD^\ast D}) - {\rm
tr}(e^{-tDD^\ast}) $
of the heat operator
is time independent, with ${\bf D} =
D+D^\ast$ acting on sections of $E\oplus F.$
   As $t\to \infty$, the supertrace
converges to the index of $D$ by a Hodge theory argument.
By parabolic regularity theory, the heat operator $e^{-t{\bf D}^\ast
{\bf D}}$ has a smooth integral kernel.   The 
pointwise supertrace of the
integral kernel of $e^{-t{\bf D}^\ast {\bf D}}$
is nonlocal, as it is
built from the eigensections of ${\bf D}^\ast {\bf D}$, which are
solutions of a global equation.  However,
by more or less standard
pseudodifferential operator techniques, the pointwise supertrace
has the asymptotic expansion 
$${\rm tr}_s e^{-t{\bf D}^\ast {\bf D}}(t,x,x) = P(x) + {\rm O}(t),$$
as $t\to 0$, 
for some local curvature invariant $P(x)= P(x,g)$, where $g$ is the
Riemannian metric.  Equating the long and short time behavior gives
${\rm Ind}\ D = \int_M P(x){\rm dvol}$.  

From the complicated construction of the asymptotics,
the curvature invariant $P(x)$ appears to contain many covariant
derivatives of the curvature, while Chern-Weil forms contain no such
derivatives.  
Either the nontrivial invariant theory
developed by Atiyah-Bott-Patodi \cite{atiyah-bott-patodi} 
and Gilkey \cite{gilkeytwo},
or the later, more direct ``supersymmetric'' approach
of Getzler \cite{getzler}  and Patodi
\cite{patodi}  shows that the covariant derivative terms
vanish and identifies $P(x) {\rm dvol} $ with the
expected Chern-Weil form: the Euler form for the Gauss-Bonnet
operator, the Todd class for the Dolbeault operator, Hirzebruch's
L-polynomial for the signature operator, and the $\hat A$-polynomial
for the Dirac operator.  
These results are sometimes called local
index theorems; proofs of various cases are in
\cite{berline-getzler-vergne}, \cite{roeI}, \cite{rosenberg-book}.
  As explained below,
the geometric  operators
all have twisted versions formed by coupling the operator to a
connection
on an auxiliary bundle.
The heat equation proof extend to the twisted case, and by
K-theory arguments the full Index Theorem reduces to
the cases of either twisted signature operators or twisted Dirac
operators.

  Bismut \cite{bismut-families} and others extended
Getzler's techniques to give similar local proofs of the Families
Index Theorem.  Here the Chern character as a differential form
constructed by Chern-Weil theory from a connection on the index bundle
is identified with the differential form version of the right hand
side of the Families Index Theorem; 
in particular, the topological pushforward is refined
to its smooth analog, integration along the fiber.   These proofs are
based on the superconnection techniques introduced in \cite{quillen2}.

Although not as popular as the heat equation approach, one can prove
the Index Theorem using the zeta functions $\zeta^{\pm}(s)$ 
of $ \Delta^+=D^\ast{ D}, \Delta^- = DD^\ast$. Here $\zeta^\pm(s) =
\sum
(\lambda_i^\pm)^{-s}$,
where $\{\lambda_i^\pm\}$ are the nonzero eigenvalues of $\Delta^\pm$
counted with multiplicity.   The smoothness of the heat kernel gives
the convergence of its trace $\sum e^{-\lambda_i^\pm t}$, which
implies a Weyl-type estimate on the growth of the $\lambda_i^\pm.$
This in turn implies the convergence of $\zeta^\pm(s)$ for ${\rm
Re}(s) >>0.$  Via a Mellin transform, the asymptotic expansion of the
heat kernel gives a meromorphic continuation of $\zeta^\pm(s)$ to all
of ${\bf C}$ with at most simple poles at prescribed integers (or
half-integers) \cite{seeley}.
In short, one obtains $\zeta^+(0) = {\rm
dim\ Ker}\ D - \int_M P^+(x)\dvol$ for some complicated local
invariant $P^+(x)$,  and similarly for $\zeta^-(0).$  It
is easy to show that $\zeta^+(0) = \zeta^-(0)$, and then supersymmetry
arguments as above yield $P^+(x) - P^-(x) = P(x)$, from which the
Index Theorem follows.

In analogy with number theoretic zeta functions, the special values or
residues of the poles of the zeta function have geometric
significance--e.g. $\zeta^\pm({\rm dim}\ M/2)$ is a multiple of the
volume of $M$.  Moreover, the zeta functions associated to locally
symmetric spaces often coincide with number theoretic zeta functions.
  However, this is as far as the analogy goes.  For example, 
it is easy to check that the zeta function for the Laplacian on the
simplest closed manifold, the circle, is a slight variation of the
Riemann zeta function, so it is not a good idea to ask a geometer if
the Riemann hypothesis holds for $\zeta^\pm(s).$   

In general, zeta functions associated to a Riemannian metric do not
have a functional equation, the main tool for studying arithmetic zeta
functions.  On the other hand, geometers can study the metric
dependence of the value of a zeta function, a technique usually
unavailable to number theorists due to rigidity results.  This
variational approach to special values of zeta functions and the
analogous L-series is a central theme in the subject.

\bigskip
\noindent \large{\bf \S2 Nonlocal invariants -- the eta invariant}

\bigskip
\normalsize

The Gauss-Bonnet theorem for a surface with boundary picks up a
correction term, the integral of the geodesic curvature over the
boundary.  Chern's generalization of Gauss-Bonnet to higher
dimensions similary has a term given by integrating a local invariant
over the boundary.  (The boundary term involves the geometric
placement of the boundary in the manifold, and so cannot be expressed
simply by a Chern-Weil-type curvature expression.)  

However, other examples of the Index Theorem do not admit such simple
generalizations.  In particular, assume the signature operator
admitted a
local correction term as above: $\sigma(M) = \int_M P_{\rm sig}(x,g) +
\int_{\partial M} Q_{\rm sig}(x,g)$.  If $\pi:M^\prime\to M$ is a 
$k$ to $1$ finite
covering,  then with respect to the pullback metric $g^\prime$ of
$g$, we would have $P_{\rm sig}(y,g^\prime) = P_{\rm sig}(\pi y,g),$
and similarly for $Q.$  
This implies $\sigma(M^\prime) = k\cdot \sigma(M)$.  However, this
formula is false in general.

Nevertheless, the dependence of the error term $\eta = \sigma(M) -
\int_M P_{\rm sig}(x,g) $ on the metric $g$ is a local invariant.  For
if $g_t$ is a curve of metrics,
%    with variation $s = (d/dt)|_{t=0} g_t$,
%   a symmetric two-tensor, 
then by Chern-Weil theory
$(d/dt)|_{t=0} \eta^{g_t}$ is a local invariant.  Of course, it is
highly
nontrivial to identify $\eta$ from this information.

In \cite{apsI}, Atiyah, Patodi and Singer gave the correction term for
the
Index Theorem for manifolds with boundary.
If we assume for simplicity that the metric 
splits as a product
on a collar of the boundary, the signature operator splits into $D =
d/dt + A$, where $t$ is the collar parameter,
$A = \pm (\ast d - d\ast)$ acting on all forms on the
boundary, and $\ast$ is the Hodge star.  The correction term depends
only on $\partial M$, so with hindsight it is reasonable to guess that
this term is built from the spectrum of the first order, selfadjoint
elliptic operator $A$.  Since $A$ has an infinite number of positive
and of negative eigenvalues, neither the heat operator nor the zeta
function of $A$ is defined.  However, in analogy with number theoretic
L-series, we can define the eta function of $A$ by
$$\eta(s) = \sum_i \frac{{\rm sgn}(\lambda_i)}{| \lambda_i|^s},$$
where the sum is over the nonzero eigenvalues of $A$.  
By a simple
scaling argument, the correction term must be $\eta(0)$ if the
correction is a spectral invariant of $A$. 

By a Mellin-type transform, it is possible to show that
$\eta(s)$ converges for ${\rm Re}(s)>>0$ and has a meromorphic
continuation to ${\bf C}$ with  $\eta(0)$ finite.  However, note that
the eta invariant does not vary smoothly with the metric, as a nonzero
eigenvalue of $A= A_t$ may pass through zero as $t$ varies, causing a
discontinuity in $\eta(s)$.  So the correction term in general must
keep track of the dimension of ${\rm Ker}\ A.$

We now state the Index Theorem for manifolds with boundary in the case
of the signature operator:

\bigskip
\noindent {\bf Theorem 1 [APS Theorem]}  {\it 
Let $M$ be a manifold with
boundary, and let $L(x,g)$ be the Chern-Weil representative of the
Hirzebruch L-polynomial  associated to a Riemannian metric $g$ on
$M$.  Assume the metric splits as a product near the boundary.  Then}
$$\sigma(M) = \int_M L(x,g) - \eta(0).$$

\bigskip

\noindent This theorem shows that the real valued {\it eta invariant}
$\eta(0)$  is nonlocal.  There are
similar results for other geometric operators.  If the metric is not a
product near the boundary, the right hand side of the last equation
also contains the integral of a local invariant over $\partial M$, as
in the Chern-Gauss-Bonnet theorem.  The text \cite{melrose}, which
contains extensive references, gives a 
proof of the APS theorem using Getzler's methods,
in which the eta invariant
appears naturally in the course of the proof.  Finally, we note that
the nonlocal nature of the eta invariant is reflected in the fact that
the APS theorem can be proved by considering a boundary value problem
for ${\bf D}^\ast {\bf D}$ with respect to nonlocal boundary
conditions. The index of $D$ is not quite equal to $\sigma(M)$, but it
involves a term which cancels with the
dimension of ${\rm Ker}\ A$, so this discontinuous term does not
appear in the APS theorem.

Needless to say, computing the eta invariant is quite difficult.
However, the fact that the variation of the eta invariant is local
leads to a computable invariant.  For let $E$ be a  hermitian bundle
over $M$ 
with  a flat unitary connection $\nabla$.  Form the coupled
signature operator $D_\nabla$ acting on $E$-valued forms by $D_\nabla
= D_{\rm sig}\otimes {\rm Id} + {\rm Id}\otimes \nabla$.  Since the
bundle is flat, $D_\nabla$ locally looks like $k$ copies of $D_{\rm
sig}$, where $k = {\rm rank}(E).$    Let $g_t$ be a curve of metrics
with variation $s = (d/dt)|_{t=0} g_t,$ a symmetric two-tensor.
Since the variations of the eta
invariants of $D_{\rm sig},D_\nabla$ are computed from the asymptotics
of the associated heat operators, it is easy to see that the
variations
$\delta_s$ of the eta invariants  satisfy
$k\cdot\delta_s
\eta_{D_{\rm sig}}(0) = \delta_s\eta_{D_\nabla}(0).$  In particular,
the difference 
$k\cdot\delta_s \eta_{D_{\rm sig}}(0) -\delta_s\eta_{D_\nabla}(0)$  is
 independent of the Riemannian metric and hence is 
a smooth invariant of the manifold $M$ and the flat bundle $E$ (or
equivalently, of the manifold $M$ and a unitary representation of
$\pi_1(M)$).  The Index Theorem for Flat Bundles
\cite{apsIII}, \cite{gilkey} gives a topological interpretation of
 this invariant.

\bigskip
\noindent \large{\bf \S3 Nonlocal invariants -- the determinant}

\bigskip
\normalsize

For geometric operators such as the Gauss-Bonnet and Dolbeault
operator, the index is the Euler characteristic of certain cohomology
groups.  If these groups vanish, the Index Theorem has nothing to say,
and secondary geometric and topological invariants appear.  This is
the case for Reidemeister torsion, which is (most easily) defined only
when the cohomology of $M$ twisted by a unitary representation $\rho$
of
$\pi_1(M)$ 
vanishes \cite{milnor}.  In particular, Reidemeister torsion is
a ``precohomological'' invariant, in that it is defined from a twisted
cochain complex; this torsion was originally introduced to distinguish
lens spaces with isomorphic cohomology rings and homotopy groups.

Ray and Singer \cite{ray-singer} noted that
Reidemeister torsion is computed from a triangulation of $M$ as a
linear combination of determinants of combinatorial Laplacians, and
conjectured that the Reidemeister torsion equals the analytic torsion,
$$T(M,g,\rho) = {1\over 2}\sum_{k=1}^n (-1)^k k \ln\ \det
\Delta_\rho^k,$$
the corresponding combinations of determinants of Laplacians on
 $k$-forms with values in the flat bundle determined by $\rho.$  

Of course, the product of the eigenvalues of these Laplacians is
infinite, so the notion of determinant must be redefined, or {\it
regularized}.  Since the determinant of an invertible finite
dimensional positive selfadjoint transformation $T$ can be computed
from its eigenvalues $\{\lambda_i\}$ by the formula $ -\ln\ \det \ T =
(d/ds)|_{s=0} \sum \lambda_i^{-s}$, Ray and Singer defined the
determinant
of a Laplacian-type operator to be $\exp(- \zeta^\prime(0)).$   In
a similar way, $\zeta(0)$ regularizes the dimension of the domain of
the Laplacian, and more importantly, the eta invariant regularizes the
number of positive minus the number of negative eigenvalues of the
operator $A$.

Based on the evidence presented in \cite{ray-singer}, Cheeger
\cite{cheeger} and M\"uller \cite{muller} proved the equality of
analytic and Reidemeister torsion.   Cheeger's proof uses surgery
techniques to reduce the problem to lens spaces, where the result was
known, while M\"uller's proof examines the convergence of
the spectral theory of the combinatorial Laplacians to that of the
smooth Laplacians as the mesh of the triangulation goes to zero.
Vishik \cite{vishik} has given  a cutting and pasting proof based on
ideas from topological quantum field theory, and 
Bismut and Zhang  \cite{bismut-zhang} have a proof based on Witten's
proof
\cite{witten-morse} of the Morse inequalities.

For a fixed Laplacian, such as the Laplacian on bundle-valued
$k$-forms, $\zeta^\prime(0)$ is a subtle spectral invariant.  In
contrast to $\zeta(0)$, it is not a local invariant, as covering space
examples show, and in contrast to $\eta(0)$, its variation is not a
local invariant.  This can be seen
by computing examples on flat tori, where the zeta function 
is a Selberg-Chowla zeta function.  
Thus there is no index theorem for flat bundles for
$\zeta^\prime(0).$  However, for operators with reasonable change
under conformal change of the metric, $\zeta^\prime(0)$ is an
invariant of the conformal class of a metric; examples include Dirac
operators and conformal Laplacians \cite{branson-orsted},
\cite{parker-rosenberg},
\cite{rosenberg}
 (cf.~\cite{ops} for a proof of the uniformization
theorem for Riemann surfaces based on the variation of
$\zeta^\prime(0)$).

In direct contrast to the original Index Theorem, which only gives
Ind $D =0$ if 
dim $M$ is odd, Reidemeister/analytic torsion is nontrivial
only if dim $M$ is odd.
Note the continued analogy with number theoretic
zeta functions, where the residue or leading order term of a special
value has significance, in that analytic torsion is
independent of the metric only in
the case where $\zeta(0)=0$.

\bigskip
\noindent\large{\bf \S4 Relations among the invariants}

\bigskip
\normalsize

All known relations among the local and nonlocal invariants involve
the
Families Index Theorem.   The first relation uncovered involves the
determinant line bundle of a family of operators, which 
is the easiest index bundle to understand topologically, but whose
geometry is nontrivial.

Motivated by questions in quantum field theory, Atiyah and Singer
\cite{as-dops}, \cite{singer}, studied the problem of canonically
defining the
determinants of a  family of invertible first order elliptic operator
$D:\Gamma(E)\to\Gamma(F)$ for different bundles $E,F.$  If $D = D_n$
were a
family of transformations, parametrized by $n\in N$,
between finite dimensional vector spaces
$V,W$, then defining the determinant of the family is equivalent to
picking a section of the line bundle 
$$\Lambda^{\rm max} V \otimes (\Lambda^{\rm max} W)^\ast \simeq
\Lambda^{\rm max} {\rm Ker}\ D \otimes (\Lambda^{\rm max}
{\rm Coker}\ D)^\ast$$
over $N$.
For the infinite dimensional vector spaces $\Gamma(E), \Gamma(F)$, 
only the right hand side of the last equation is defined, so defining
the determinant of the family reduces to
trivializing the determinant line bundle ${\sc Det}\ D = \Lambda^{\rm
max}{\sc Ker }\ D 
\otimes (\Lambda^{\rm max}{\sc Coker}\ D)^\ast$ of the index bundle.
Thus there are topological obstructions to the solution given by the
Families Index Theorem.
However, even if the  bundle is trivial, there is no canonical
trivialization, unless the bundle comes with a canonical nonvanishing
section, or
equivalently with a canonical flat connection with vanishing
holonomy.  

In \cite{quillen}, Quillen introduced a canonical hermitian
metric on the determinant line bundle for the family  of twisted
Dolbeault
operators over a Riemann surface.  In brief, {\sc Det} $D$ has a
canonical
section det $D$ which is locally given by
$\det \ D = v^1\wedge\ldots\wedge v^k\otimes
(Dv^1\wedge\otimes\ldots\wedge
Dv^k)^\ast$,
where $\{v^i\}$ is an
orthonormal basis of eigensections of $D^\ast D$ with eigenvalues
$\lambda_i$ 
lying below some value $a$.
%Since the
%dimension of the kernel and cokernel may jump discontinuously, the
%only way to smoothly assign a length to this section is to set
The Quillen metric on {\sc Det} $D$ is defined by
$$||\det\ D||^2 = (\prod \lambda_i)e^{-\zeta_{D^\ast
D,a}^\prime(0)},$$
where the zeta function $\zeta_{D^\ast D,a}(s)$ is built from the
eigenvalues
above $a$.
This formula is independent of the choice of $a$, and is
motivated by the easy equality
$||\det\ D||^2 = \det\ D^\ast D$
for $D$ a transformation of finite dimensional vector spaces with
$||D^\ast D||< a.$ 
In Quillen's holomorphic example, {\sc Det}
$D$ has a canonical connection, whose curvature is nonzero in general;
thus there is a geometric obstruction to defining the determinant of
the family of operators, even if the topological obstruction (the
first Chern class of {\sc Det} $D$) vanishes.

The holomorphic category is too rigid to handle the general case of
twisted Dirac operators.  
In \cite{bismut-freedI}, Bismut and Freed defined an analogous
connection on the determinant line bundle for a family
of twisted
Dirac operators, and proved that the curvature of this connection,
which (i) is a combination of second variations of $\zeta^\prime(0)$,
and (ii) 
represents the first Chern class, is
precisely the degree two component of 
the differential form in the
local Families Index Theorem.
In light of the nonlocal nature of the first variation of
$\zeta^\prime(0)$, this localization of the curvature is especially
striking.

Computing the holonomy of the Quillen metric is more difficult, as
this involves integrating the connection one-form around loops in the
base manifold. This is a nonlocal computation on a general bundle, 
and is all the more nonlocal here, since the connection one-form is
essentially the first variation of 
$\zeta^\prime(0)$.  Of course, if the family of
operators is constant, the holonomy is trivial.  In
\cite{witten}, Witten gave a heuristic argument for the limit of the
holonomy as the metric on the base manifold is scaled up; this
so-called adiabatic limit 
corresponds to taking a limit as the operators get closer and closer
to being constant.

The geometric framework for Witten's argument 
consists of
a  Riemannian submersion $\pi:X\to N$ 
 of Riemannian manifolds with fiber
diffeomorphic to the even dimensional closed  manifold 
$M$.  Assuming $X,M,N$ are spin,
we have a family of Dirac operators $D_n, n\in N,$ acting on positive
spinors on
$\pi^{-1}(n).$  (For simplicity, we will not discuss the twisted
case.) The
holonomy  around a loop $\gamma \subset N$ 
is formally controlled by the phase of  the determinant
of the associated Dirac operator {\bf D}
on all spinors on $\pi^{-1}(\gamma).$ 
This phase is formally $(-1)$ to the number
 of negative
eigenvalues of {\bf D}, but this number
must of course be regularized. If {\bf
D} were a finite dimensional transformation, the number of negative
eigenvalues would be
exactly $(\zeta_{\bf D}(0) - \eta_{\bf D}(0))/2,$
 so in general $\exp(i\pi\eta(0))$ enters in the
the phase of the determinant.  This informal argument was
presented rigorously in \cite{bismut-freedII}, \cite{dai-freed}.

\bigskip
\noindent {\bf Theorem 2 }  {\it 
 Let $\{D\}= \{D_n\}$ be a family of Dirac operators,
parametrized by $n$ in a compact Riemannian manifold $N$ as above.
The curvature of the determinant line bundle ${\sc Det} \ D$ is given
by the degree two component of 
$ \int_M \hat A(D_n,g)$, where $\hat A(D_n,g)$ is the $\hat
A$-polynomial as
a differential form associated
to the family $\{D\}$ and a metric $g$ on $N$, and
$\int_M $ denotes integration over the fiber.  The holonomy of a loop
 $\gamma\subset N$ is $(-1)^{{\rm Ind}\ {\bf D}} \exp(i\pi{\bf
\eta})$,
 where ${\bf \eta}$ is the adiabatic limit of $(\eta_{\bf D}(0) + {\rm
 dim\ Ker} \ {\bf D})/2.$}

\bigskip
\noindent In summary, in the determinant line bundle picture, the
 integral of
the first variation of
$\zeta^\prime(0)$ around a loop  relates the holonomy of the
Quillen-Bismut-Freed connection to
the eta invariant of a 
Dirac operator over the loop, 
while the Families Index
Theorem computes the curvature, which depends on the second variation
of $\zeta^\prime(0)$.

Certainly there is a loss of information in passing from the index
bundle to the determinant line bundle.  
The relation between the eta invariant and the Families Index Theorem
at the level of the index bundle
has been made more precise in work of Bismut and Cheeger
\cite{bismut-cheegerI}, \cite{bismut-cheegerII}.
For this theory
let $\pi:X\to N$ be a Riemannian submersion 
with fiber
diffeomorphic to the even dimensional  manifold with boundary
$M$.   We assume that $X,N,M$ are spin, and that the kernel of the
Dirac operator $D = D_n$ restricted to the boundary
of each fiber has constant rank.  Then the local Families Index
Theorem for
Manifolds with Boundary is the equality of differential forms
\[{\rm ch}({\sc Ind} \ D) = \int_M\hat A(D_n,g) -\tilde \eta,\]
where the {\it eta form}
$\tilde\eta$ is an even dimensional form on $N$ canonically
constructed from $D|_{\partial X},$ and $g$ is the metric on $X$.
  (There is a similar result if $M$
is odd dimensional.)
In particular, the degree zero component of $\tilde \eta$ satisfies
$\tilde\eta_0 (n)=  \eta(0)$, the eta invariant of $D_n$, and as such
this Index Theorem is a strong generalization of the APS theorem.

In addition, $\tilde \eta$ satisfies the transgression formula on
$\partial X$
\[d\tilde\eta = {\rm ch}({\sc ind}\ D_{\partial M}) -\int_{\partial
M} \hat A(D_{\partial X}, g|_{\partial X})\]
when $M$ is odd dimensional and the dimension of ${\rm Ker}\
D_{\partial M}$ is constant.  Since the right hand side is a local
invariant, the transgression formula generalizes the result that the
variation of the eta invariant is local.  The corresponding formula
for
$M$ even dimensional simply reads $d\tilde\eta =0,$ so $\tilde\eta$
determines an even dimensional cohomology class in $N$ associated
directly to $\partial X$.  

In recent work \cite{bismut-lott},  
Bismut and Lott have proved an analogous theorem which makes precise
the sense in
which analytic torsion appears as a secondary invariant once both
sides of the index theorem vanish {\it a priori}.
Let $E$ be a flat bundle over $X$, where as
above we have the Riemannian submersion $\pi:X\to N$ with fiber $M$.
For $k$ odd, 
let $c_k(E)\in H^k(X; {\bf R})$ denote the Kamber-Tondeur classes
of $E$ \cite{kamber-tondeur}; these classes have a Chern-Weil
description given below.
Let $e(TM/N)$ be the Euler
class of $TM/N$, the tangent bundle along the fibers. 
Let $H^p(M;E|_M)$ denote the flat vector bundle over $N$ whose fiber
at $n\in N$ is the cohomology group $H^p(M_n;E|_{M_n})$.  
Then we have the following cohomological result.

\bigskip
\noindent {\bf Theorem 3 }  {\it For any odd integer k,
\[ \pi_\ast (e(TM/N) c_k(E)) = \sum_{p=0}^{{\rm
dim}\ M} (-1)^p c_k(H^p(M;E|_M)) \]
as classes in $H^k(N;{\bf R}).$ }
 
\bigskip
\noindent This theorem is a $C^\infty$ (global)
analogue of a (local) version of
the Grothendieck-Riemann-Roch theorem for holomorphic
submersions
due to Bis\-mut, Gillet and Soul\'e
\cite{bismut-gillet-soule}, and as such is a families-type index
theorem.  

The striking feature of the proof of  Theorem 3 is its 
local, geometric nature; a
global, topological version appeared subsequently \cite{dww}.
To state the more precise geometric version, pick a hermitian metric
$h$ on $E$.  Since $E$ is flat with respect to a connection
$\nabla^E$, 
locally there exists a basis $\{e_i\}$ of the fibers of $E$
with $\nabla^E e_i = 0.$  With respect to this basis, $h$ is a
function on $X$ with values in Hermitian matrices.  Define an
End($E$)-valued one-form on $X$ by  $\omega(E,h) = h^{-1}dh.$ For odd 
positive integers $k$, set 
$$c_k(E,h) = (2\pi i)^{-(k-1)/2}2^{-k}{\rm Tr}(\omega^k(E,h)),$$
where $\omega^k$ is $\omega$ wedged with itself $k$ times.
The $k$-form $c_k(E,h)$ is closed and its cohomology class equals
$c_k(E).$
Let 
$\Omega$ be the curvature of
a unitary connection $\nabla$
on $TM/N$, and let $e(\Omega)$ be the Chern-Weil representative of
$e(TM/N).$ 
The fibers
$H^p(M_n;E|_{M_n})$ 
inherit an inner product $h^{H^p}$
from the $L^2$ inner
product on $\Gamma(TM/N_n\otimes E|_{M_n})$, and we denote the
associated 
Chern-Weil representative of $c_k(H^p(M;E|_M))$
by 
$ c_k(H^p(M;E|_M),h^{H^p})$.
Recall that $g$ is the metric on $X$.

\bigskip
\noindent {\bf Theorem 4}
{\it There exist $(k-1)$-forms 
\[{\cal T}_{k-1}= {\cal T}_{k-1}
(g,h)\]
 on $N$ such that

\noindent (i)  for any odd integer k,}
\[
d{\cal T}_{k-1} = 
\int_M e(\Omega)c_k(E,h) 
-\sum_{p=0}^{{\rm
dim}\ M}(-1)^p c_k(H^p(M;E|_M),h^{H^p})\]
{ \it as differential k-forms;

\noindent (ii)  if the cohomology groups $H^p(M_n;E|_{M_n})$ vanish
for all $p$
 and all $n\in N$, and if $M$ is odd dimensional, then the forms
 ${\cal
T}_{k-1}$ are closed, and the class $[{\cal T}_{k-1}]\in
 H^{k-1}(N;{\bf
R})$ is independent of the choices of $g,h.$  In
particular, $[{\cal T}_0]$ is (represented by) the locally constant
function which is half the analytic/Reidemeister torsion of the pair}
$(M_n, E|_{M_n}).$ 

\bigskip

This result has several striking features.  First, it has Theorem 3 as
an immediate corollary.  Second,
it realizes the
an\-a\-lyt\-ic/Reid\-e\-mei\-ster 
torsion as the degree zero component of a sophisticated local
familes-type index theorem.  Third, 
note that the first statement in (ii) follows immediately from (i),
since dim $M$ odd implies 
$\int_M e(\Omega)c_k(E,h)$ is an even form and so has no
component in degree $k$.  Thus, as expected, the so-called  {\it 
higher  torsion
forms} ${\cal T}_{k-1}$ have cohomological
significance only when the
cohomological information in the index-type  Theorem 3 is
trivial.   

Of course, one expects that the cohomology classes of the higher
torsion forms  should equal some
combinatorially defined ``higher Reidemeister torsion.''  Such a
torsion has been defined by Igusa and Klein \cite{i}, \cite{k}, and
the two higher torsions agree
in the few
cases where they have been computed.
At present, however, the equality in full generality
 is unknown.  

In a related direction, Lott \cite{lott} has used an extension
$\widetilde{\cal T}$
of the higher torsion forms to study $\pi_i({\rm
Diff}(Z))\otimes {\bf Q}$, where Diff$(Z)$ denotes the diffeomorphism
group of a
closed smooth $K(\Gamma,1)$ manifold $Z$.  
For a given $\alpha:S^i\to {\rm Diff}(Z)$, the standard clutching
construction of  gluing two $(i+1)$-balls times $Z$ along their
boundaries via $\alpha$  produces
 a manifold $M_\alpha$ fibering over $S^{i+1}$ with fiber $Z$.   
Under certain technical hypotheses, $\int_{S^{i+1}} \widetilde{\cal
T}(M_\alpha)$ lies in $H_*(\Gamma;{\bf C}).$  In fact, there is an
explicit
conjectured formula for $\pi_i({\rm
Diff}(Z))\otimes {\bf Q}$ 
in terms of
$H_k(\Gamma,{\bf C})$.  Lott's map $\int_{S^{i+1}}\widetilde{\cal
T}(M_\cdot):  
\pi_i({\rm
Diff}(Z))\to H_\ast(\Gamma;{\bf C})$
 gives the first analytic
evidence for this old conjecture (which is known in several cases),
although at   
this point it is not known if this map is an isomorphism.
In contrast, 
 the homeomorphism group of $Z$ has vanishing  rational homotopy
groups, 
so  it seems natural to approach the conjecture via a smooth,
i.e.~analytic, construction.

In summary, strong relations between the nonlocal invariants and
classical index theory are now firmly established.  These relations
have applications to algebraic K-theory indicated in
\cite{bismut-lott}, \cite{dww}, \cite{lott},
as well as to the subject of
index theory on open manifolds (e.g.~\cite{ads},
\cite{bismut-cheegerIII}, \cite{burg-etc}, \cite{mathai-carey},
\cite{mullerII},
\cite{roeII}),
that go beyond the limits of this article.

\bibliographystyle{amsplain}
\bibliography{art}

\end{document}